\newif\ifpdf
\ifx\pdfoutput\undefined
\pdffalse % we are not running PDFLaTeX
\else
\pdfoutput=1 % we are running PDFLaTeX
\pdftrue \fi

\newif\iffinal
\finalfalse	% Not final version
\finaltrue % Final version

\documentclass[reqno,twoside,11pt]{amsart}
\iffinal\else\usepackage[notref,notcite]{showkeys}\fi
\usepackage{amsmath}
\usepackage{amsfonts}
\usepackage{amssymb}
\usepackage{verbatim}
\usepackage{epsfig}
\IfFileExists{epsf.def}{\input epsf.def}{\usepackage{epsf}}

\IfFileExists{myowntimes.sty}{\usepackage{myowntimes}\usepackage{mathrsfs}}
	{\usepackage{times}\usepackage{mathrsfs}}
\DeclareFontFamily{OT1}{eusb}{} \DeclareFontShape{OT1}{eusb}{m}{n} {<5> <6> <7> <8> <9> <10> <11> <12> <14.4> eusb10}{}
\DeclareMathAlphabet{\eusb}{OT1}{eusb}{m}{n}

\DeclareFontFamily{OT1}{eusm}{} \DeclareFontShape{OT1}{eusm}{m}{n} {<5> <6> <7> <8> <9> <10> <11> <12> <14.4> eusm10}{}
\DeclareMathAlphabet{\eusm}{OT1}{eusm}{m}{n}

\DeclareFontFamily{OT1}{eufm}{} \DeclareFontShape{OT1}{eufm}{m}{n} {<5> <6> <7> <8> <9> <10> <11> <12> <14.4> eufm10}{}
\DeclareMathAlphabet{\mathfrak}{OT1}{eufm}{m}{n}

\DeclareFontFamily{OT1}{fraktura}{}
\DeclareFontShape{OT1}{fraktura}{m}{n} {<5> <6> <7> <8> <9> <10> <11> <12> <13> <14.4> [1.1] eufm10}{}
\DeclareMathAlphabet{\fraktura}{OT1}{fraktura}{m}{n}

\DeclareFontFamily{OT1}{cmfi}{} \DeclareFontShape{OT1}{cmfi}{m}{n} {<5> <6> <7> <8> <9> <10> <11> <12> <13> <14.4> [0.9] cmfi10}{}
\DeclareMathAlphabet{\cmfi}{OT1}{cmfi}{b}{n}

\DeclareFontFamily{OT1}{cmss}{} \DeclareFontShape{OT1}{cmss}{m}{n} {<5> <6> <7> <8> <9> <10> <11> <12> <13> <14.4> cmss10}{}
\DeclareMathAlphabet{\cmss}{OT1}{cmss}{m}{n}
\newtheoremstyle{theorem}{1.5ex}{1.5ex}{\itshape\rmfamily}{} {\bfseries\rmfamily}{}{2ex}{}

\newtheoremstyle{def}{1.5ex}{1.5ex}{\slshape\rmfamily}{} {\bfseries\rmfamily}{}{2ex}{}

\newtheoremstyle{rem}{1.3ex}{1.3ex}{\rmfamily}{} {\itshape}
{} {1.5ex}{}

\newenvironment{proofsect}[1] {\vskip0.1cm\noindent{\rmfamily\itshape#1.}}{\qed\vspace{0.15cm}}%{\newline\vspace{0.15cm}}

\theoremstyle{theorem}
\newtheorem{theorem}{Theorem}[section]
\newtheorem{lemma}[theorem]{Lemma}
\newtheorem{proposition}[theorem]{Proposition}

\newtheorem*{Main Theorem}{Main Theorem.}
\newtheorem{corollary}[theorem]{Corollary}
\newtheorem*{special theorem}{Lindeberg-Feller Theorem for Martingales}

\theoremstyle{def}
\newtheorem{definition}[theorem]{Definition}

\theoremstyle{rem}

\newtheorem{remarks}[theorem]{{\itshape Remarks}}

\numberwithin{equation}{section}

\renewcommand{\section}{\secdef\sct\sect}
\newcommand{\sct}[2][default]{\refstepcounter{section}
\addcontentsline{toc}{section}
{{\tocsection {}{\thesection}{\!\!\!\!#1\dotfill}}{}}
\vspace{0.7cm}
\centerline{ %\large
\scshape\arabic{section}.\ #1} \nopagebreak \vspace{0.2cm}}
\newcommand{\sect}[1]{
\vspace{0.4cm} \centerline{\large\scshape\rmfamily #1}
\vspace{0.2cm}}

\renewcommand{\subsection}{\secdef\subsct\sbsect}
\newcommand{\subsct}[2][default]{\refstepcounter{subsection}
\addcontentsline{toc}{subsection}
{{\tocsection{\!\!}{\hspace{1.2em}\thesubsection}{\!\!\!\!#1\dotfill}}{}}
\nopagebreak\vspace{0.45\baselineskip} {\flushleft\bf
\arabic{section}.\arabic{subsection}~\bf #1.~}
\\*[3mm]\noindent
\nopagebreak}
\newcommand{\sbsect}[1]{\vspace{0.1cm}\noindent
\textbf{#1.~}\vspace{0.1cm}}

\renewcommand{\subsubsection}{%
\secdef \subsubsect\sbsbsect}
\newcommand{\subsubsect}[2][default]{%
\refstepcounter{subsubsection} 
\addcontentsline{toc}{subsubsection}{{\tocsection{\!\!}
{\hspace{3.05em}\thesubsubsection}{\!\!\!\!#1\dotfill}}{}}
\nopagebreak
\vspace{0.15\baselineskip} \nopagebreak {\flushleft\rmfamily
\itshape\arabic{section}.\arabic{subsection}.\arabic{subsubsection}
\ \rmfamily #1\/.}\ }
\newcommand{\sbsbsect}[1]{\vspace{0.1cm}\noindent
\rmfamily \itshape
\arabic{section}.\arabic{subsection}.\arabic{subsubsection} \
\sffamily #1\/.\ }

\iffinal

\else

\fi

\renewcommand{\caption}[1]{%
\vglue0.5cm
\refstepcounter{figure}
\begin{minipage}{0.9\textwidth}\small {\sc Figure~\thefigure. }#1\end{minipage}}

\newcommand{\textd}{\text{\rm d}\mkern0.5mu}
\newcommand{\texti}{\text{\rm  i}\mkern0.7mu}
\newcommand{\texte}{\text{\rm e}}

\newcommand{\C}{\mathbb C}

\newcommand{\E}{\mathbb E}

\newcommand{\scrE}{\mathscr{E}}

\newcommand{\scrG}{\mathscr{G}}

\newcommand{\scrV}{\mathscr{V}}
\newcommand{\scrW}{\mathscr{W} }

\newcommand{\twoeqref}[2]{(\ref{#1}--\ref{#2})}
\newcommand{\cc}{{\text{\rm c}}}

\newcommand{\alphac}{\alpha_\cc}
\def\myffrac#1#2 in #3{\raise 2.6pt\hbox{$#3 #1$}\mkern-1.5mu\raise 0.8pt\hbox{$#3/$}\mkern-1.1mu\lower 1.5pt\hbox{$#3 #2$}}
\newcommand{\ffrac}[2]{\mathchoice%
	{\myffrac{#1}{#2} in \scriptstyle}
	{\myffrac{#1}{#2} in \scriptstyle}
	{\myffrac{#1}{#2} in \scriptscriptstyle}
	{\myffrac{#1}{#2} in \scriptscriptstyle}
}

\newcommand{\RE}{\mkern2mu\text{\rm Re}\mkern1mu}

\providecommand{\p}{\ffrac{\alpha}{n}}

\title[Large-deviations and complete-graph percolation]
{\Large Large-deviations/thermodynamic approach\\to percolation on the complete graph}
\author[M.~Biskup, L.~Chayes and S.A.~Smith]{Marek Biskup,\, Lincoln Chayes\, and\, S.~Alex Smith}

\begin{document}
\thanks{\hglue-4.5mm\fontsize{9.6}{9.6}\selectfont\copyright\,2006 by M.~Biskup, L.~Chayes, S.A.~Smith. Reproduction, by any means, of the entire
article for non-commercial purposes is permitted without charge.\vspace{2mm}}
\maketitle

\vspace{-5mm}
\centerline{\textit{Department of Mathematics, University of California at Los Angeles}}

\begin{abstract}
We present a large-deviations/thermodynamic approach to the classic problem of percolation on the complete graph. Specifically, we determine the large-deviation rate function for the probability that the giant component occupies a fixed fraction of the graph while all other components are ``small.'' One consequence is an immediate derivation of the ``cavity'' formula for the fraction of vertices in the giant component. 
As a by-product of our analysis we compute the large-deviation rate functions for the probability of the event that the random graph is connected, the event that it contains no cycles and the event that it contains only ``small'' components.
\end{abstract}

\section{Introduction}
\noindent
For physical systems, mean-field theory often provides a qualitatively correct description of ``realistic behavior.'' The corresponding analysis usually begins with the derivation of so called mean-field equations which are self-consistent relations involving the physical quantity of primary interest and the various parameters of the model. This approach may be realized and, to some extent, justified mathematically by considering the model on the complete graph where each constituent interacts with all others. 

As an example, let us consider the Ising model on a complete graph~$K_n$ of~$n$ vertices. Here we have a collection of~$\pm1$-valued random variables~$(\sigma_i)_{i=1}^n$ which are distributed according to the probability measure~$\mu_n(\{\sigma\})=\texte^{-\beta H_n(\sigma)}/Z_{n,\beta}$, where
\begin{equation}
H_n(\sigma)=-\frac1n\sum_{i,j=1}^n\sigma_i\sigma_j-h\sum_{i=1}^n\sigma_i
\end{equation}
and where~$\beta,h$ are parameters. The relevant physical quantity is the \emph{empirical magnetization}, $m_n(\sigma)=n^{-1}\sum_{i=1}^n\sigma_i$. In terms of this quantity, $H_n(\sigma)=-\frac12 n[m_n(\sigma)]^2-h m_n(\sigma)$ and so
\begin{equation}
\E_n(\sigma_1|\sigma_j\colon j\ne1)=\tanh\bigl[\beta(m_n(\sigma)+h)\bigr]+O(\ffrac1n).
\end{equation}
This permits the following ``cavity argument:'' Supposing that~$m_n$ tends, as~$n\to\infty$, to a value~$m_\star$ in probability, we have that $m_\star=\lim_{n\to\infty}\E_n(\sigma_1)$ obeys
\begin{equation}
\label{m-tanh-m}
m_\star=\tanh\bigl[\beta (m_\star+h)\bigr].
\end{equation}
This is the \emph{mean-field equation} for the (empirical) magnetization. Of course, the concentration of the law of~$m_n$ still needs to be justified; cf~\cite{Ellis} for details.

In the context of percolation~\cite{Grimmett}, the relevant mean-field model goes under the name the Erd\"os-Renyi Random Graph. Here each edge of~$K_n$ is independently occupied with probability~$\p$, where~$0\le \alpha<\infty$, and vacant with probability~$1-\ffrac\alpha n$. The relevant ``physical'' quantity is the \emph{giant-component density}~$\varrho_\star$, i.e., the limiting fraction of the vertices that belong to the giant component of the graph. The corresponding mean-field equation,
\begin{equation}
\label{1.1}
\varrho_\star=1-\texte^{-\alpha \varrho_\star},
\end{equation}
is also readily derived from heuristic ``cavity'' considerations.   As is well known, $\varrho_\star=0$ is the only solution for~$\alpha\le\alphac=1$, while for~$\alpha>\alphac$ there is another, strictly positive solution. This solution tends to zero as~$\alpha\downarrow\alphac$; hence we may speak of a continuous transition.

While \twoeqref{m-tanh-m}{1.1} are indeed straightforward to derive, matters at the level of mean-field equations are not always satisfactory; the problem being the existence multiple solutions. As it turns out, for the percolation model (as well as the $k$-core percolation) the proper choice is always the \emph{maximal} solution, but prescriptions of this sort generically fail, e.g., for the Ising model \eqref{m-tanh-m} with~$h<0$ and, as often as not, whenever there is a first-order transition. Thus, one is in need of an additional principle which determines which of the solutions is relevant.

The existing mathematical approach to these difficulties---e.g., for percolation~\cite{Erdos-Renyi}, see also~\cite{Alon-Spencer,Bollobas,Janson-Luczak-Rucinski}, or the $k$-core~\cite{Pittel-Spencer-Wormald}---is to work with sufficient precision until the mean-field conclusions are rigorously established. Another approach---which admits some prospects of extendability beyond the complete graph~\cite{BC,BCN}---is to supplement the picture by the introduction of the \emph{mean-field free-energy function}. For the Ising model, this is a function $m\mapsto\Phi_{\beta,h}(m)$ such that
\begin{equation}
\mu_n\bigl(m_n(\sigma)\approx m\bigr)=\texte^{-n\Phi_{\beta,h}(m)+o(n)},\qquad n\to\infty,
\end{equation}
i.e., $m\mapsto\Phi_{\beta,h}(m)$ is the large-deviation rate function for the probability of observing the event $\{m_n(\sigma)\approx m\}$. 
This spells the end of the story from the perspective of probability and/or theoretical physics: One seeks the minimum of the free energy function, setting its derivative to zero yields the mean-field equations with the irrelevant solutions corresponding to the local extrema which are not absolute minima; see again~\cite{Ellis}.

\begin{figure}[t]
{\epsfig{figure=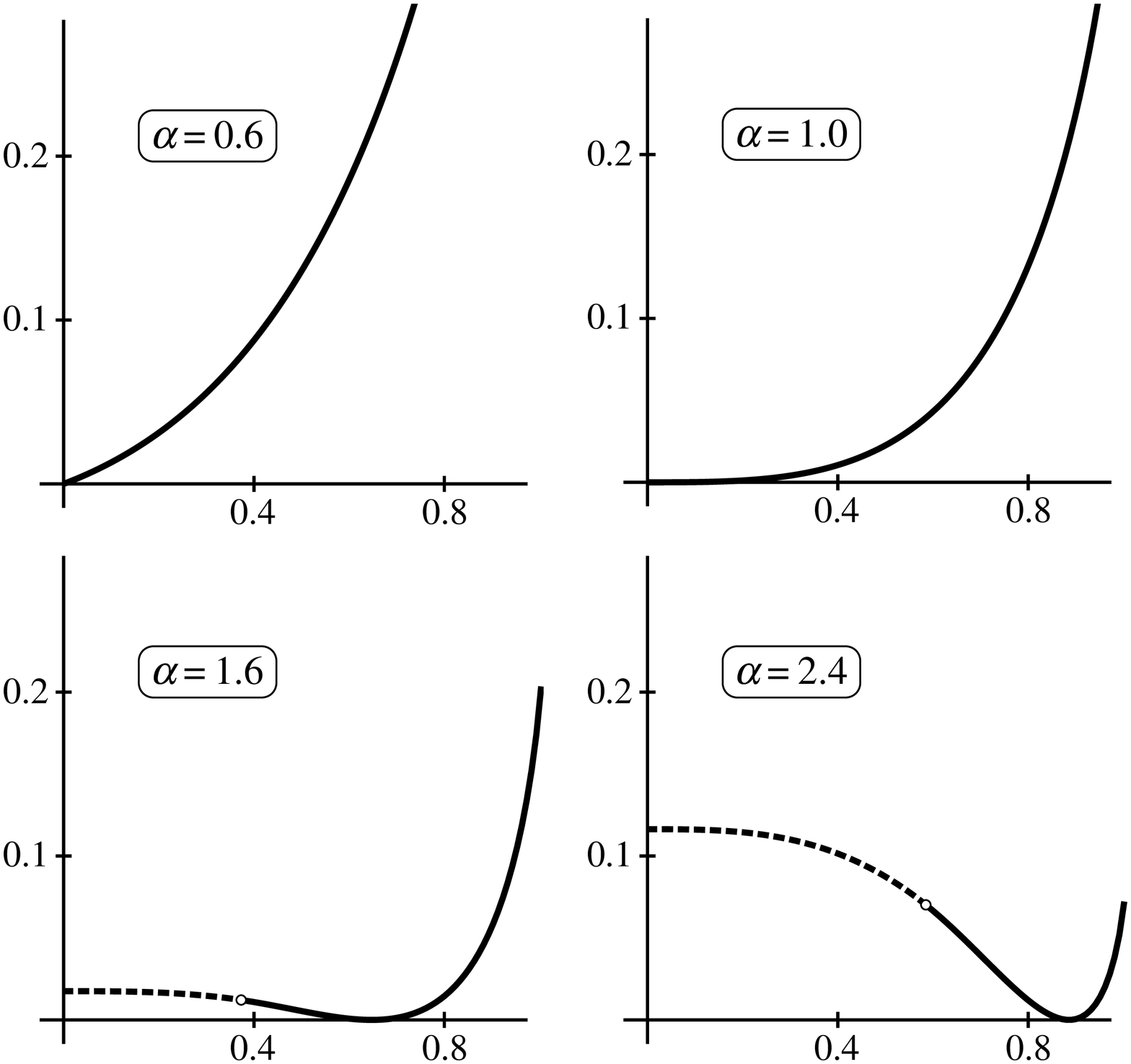, width=3.7in}}
\caption{
The graph of the free energy function $\varrho\mapsto\Phi(\varrho,\alpha)$ for four distinct values of~$\alpha$. For~$\alpha\le1$, the function is minimized by~$\varrho=0$, while for~$\alpha>1$ the unique minimum occurs at~$\varrho>0$. In any case, the minimizer is the maximal positive solution for~$\varrho_\star$ from~\eqref{1.1}. The dashed portion of the graphs for~$\alpha=1.6$ and~$2.4$ marks the part where the background contribution,~$\Psi(\alpha(1-\varrho))$, to $\Phi(\varrho,\alpha)$ in~\eqref{2.4} is strictly positive. This rules out the zero solution to \eqref{1.1} for all~$\alpha>1$.
}
\label{fig1}
\end{figure}  

The free-energy approach to mean-field problems has met with success in Ising systems and, to some extent, it has been applied to the Potts and random-cluster models~\cite{Bollobas-Grimmett,Monasson-at-al,Costeniuc-at-al,Luczak-Luczak}.
However, no attempt seems to have been made to extend this technology to ``purely geometrical'' problems on the complete graph, specifically, ordinary percolation or~$k$-core percolation.
The purpose of this note is to derive the large-deviation rate function for the event that the random graph contains a fraction~$\varrho$ of vertices in ``large'' components. As we will see, the function has a unique minimum for all~$\alpha$ which coincides with the ``correct'' solution of \eqref{1.1}.  We do not necessarily claim that the resultant justification of this equation is easier than which already exists in the literature.  However, the picture presented here provides some additional insights into the model while the overall approach indeed admits the possibility of generalizations.

\section{Main results}
\noindent
Consider the set of vertices~$\scrV=\{1,\dots,n\}$ and let~$(\omega_{kl})_{1\le k<l\le n}$ be a collection of i.i.d.~random variables taking value one with probability~$p$ and zero with probability~$1-p$. Let~$\scrE=\scrE(\omega)$ be the (random) set~$\{(k,l)\colon 1\le k<l\le n,\,\omega_{kl}=1\}$. In accord with the standard notation, cf~\cite{Alon-Spencer,Bollobas}, we will use~$\scrG(n,p)$ to denote the undirected graph with vertices~$\scrV$ and edges~$\scrE$. Of particular interest are the cases where~$p$ decays to zero proportionally to~$\ffrac1n$. Since these are the only problems we will consider, let us set, for once and all, $p=\p$ for some fixed~$\alpha\in[0,\infty)$. We will denote the requisite probability measure by~$P_{n,\alpha}$.

In order to state our main theorems, we need to introduce some notation. First, consider the standard entropy function
\begin{equation}
%\label{}
S(\varrho)=\varrho\log\varrho+(1-\varrho)\log(1-\varrho)
\end{equation}
and let
\begin{equation}
%\label{}
\pi_1(\alpha)=1-\texte^{-\alpha}.
\end{equation}
In addition, consider the function
\begin{equation}
\label{eq-Psi}
\Psi(\alpha)=\biggl(\log\alpha-\frac12\Bigl[\alpha-\frac1\alpha\Bigr]\biggr)\wedge 0
\end{equation}
and note that~$\Psi(\alpha)<0$ if and only if~$\alpha>1$. Finally, let us also define
\begin{multline}
\label{2.4}
\qquad
\Phi(\varrho,\alpha)=S(\varrho)-\varrho\log\pi_1(\alpha\varrho)
\\
-(1-\varrho)\log\bigl[1-\pi_1(\alpha\varrho)\bigr]
-(1-\varrho)\Psi\bigl(\alpha(1-\varrho)\bigr).
\qquad
\end{multline}
Then we have:

\begin{theorem}
\label{thm-main}
Consider $\scrG(n,\ffrac\alpha n)$ and let~$\scrV_r$ be the set of vertices that are in connected components of size larger than~$r$. 
Then for every~$\varrho\in[0,1]$,
\begin{equation}
%\label{}
\lim_{\epsilon\downarrow0}
\lim_{n\to\infty}P_{n,\alpha}\bigl(|\scrV_{\epsilon n}|=\lfloor\varrho n\rfloor\bigr)^{1/n}
=\texte^{-\Phi(\varrho,\alpha)}.
\end{equation}
\end{theorem}

An inspection of Lemma~\ref{cor6.2} reveals that, conditional on $\{|\scrV_{\epsilon n}|=\lfloor\varrho n\rfloor\}$, with~$\epsilon>0$, there will be only one ``large'' component with probability tending to one as~$n\to\infty$.
 
\smallskip
Fig.~1 shows the graph of~$\Phi$ for various values of~$\alpha$ which is archetypal of free-energy functions in complete graph setting. The figure indicates a unique global minimum; direct, albeit arduous differentiation of \eqref{2.4} yields the fact that all local extrema satisfy the mean-field equation~\eqref{1.1}. The extremum at~$\varrho=0$ is ruled out for~$\alpha>1$ by noting that, under these conditions, the last term in \eqref{2.4} is strictly positive.

The corresponding conclusion may also be extracted from the following probabilistic argument: Let~$m=\lfloor \varrho n\rfloor$ and note that~$e^{-nS(\varrho)}$ is then the exponential growth-rate of~$\binom nm$. This allows us to write
\begin{equation}
\label{2.6}
\texte^{-n\Phi(\varrho,\alpha)}=e^{o(n)}
\binom nm\bigl[\pi_1(\alpha\varrho)\bigr]^m
\bigl[1-\pi_1(\alpha\varrho)\bigr]^{n-m}
\,\texte^{(n-m)\Psi(\alpha(1-\varrho))}.
\end{equation}
Neglecting the~$\Psi$-term (which provides a lower bound on~$\Phi$), one sees a quantity reminiscent of binomial distribution. Well known results on the latter inform us that the right-hand side is exponentially small unless
\begin{equation}
%\label{}
\pi_1(\alpha\varrho)\approx\frac mn,
\end{equation}
i.e., unless~$\varrho$ satisfies the mean-field equation~\eqref{1.1}.
If~$\Psi$ is set to zero, there are degenerate minima for~$\alpha>1$; however, the $\Psi$-function will lift the degeneracy and, in fact, create a local \emph{maximum} at~$\varrho=0$ once~$\alpha>1$. Meanwhile, in the region of the maximal solution, $\Psi$ has vanished and the above mentioned approximation is exact. 

\begin{remarks}
(1) A closely-related, but different problem to the one treated above has previously been studied using large-deviation techniques. Indeed, in~\cite{OConnell}, O'Connell derived the large-deviation rate function for the event that the \emph{largest} connected component is of size about~$\kappa n$. Note, however, that this does not restrict the total volume occupied by these component. For~$\kappa$ close to~$\varrho_\star$ from \eqref{1.1}---explicitly, as long as the complement of the large component has effective~$\alpha$ less than~$1$---O'Connell's rate function coincides with ours. But once~$\kappa$ is sufficiently small, his conditioning will lead to the creation of several large components whose total volume is such that their complement is effectively subcritical.  Consequently, O'Connell never needs to address the central issue of our proof; namely, the decay rate of the probability that supercritical percolation has no giant components. (This is what gives rise to the term~$\Psi$ in \eqref{2.4} and the dashed portion of the graph in Fig.~\ref{fig1}.) In fact, his rate function is basically a concatenation of many scaled copies of the undashed portion of the graph in Fig.~\ref{fig1}.

(2) While the $\Psi$-term in \eqref{2.4} has a non-trivial effect on the large-deviation questions studied here, it does not play any role for events whose probability is of order unity (or is subexponential in~$n$). This is because $\Psi$ ``kicks in'' only for $\varrho$ away from the minimizing value. This is not the case for the $k$-core where the corresponding large-deviation analysis~\cite{BCS-kcore} suggests that the analogous term ``kicks in'' right at the minimizer and may even affect the fluctuation scales. One way to bring~$\Psi$ out of the ``realm of exponentially-improbable'' for percolation would be to give each configuration a weight suppressing large components. However, we will not pursue these matters in the present note.

(3)
Our control of the rate function is not sharp enough to provide a detailed description of the critical region, i.e., the situations when
$\alpha=1+O(n^{-1/3})$. The corresponding analysis of the scaling phenomena inside the ``critical window'' has been performed in~\cite{Bollobas-crit,Luczak,Luczak-at-al,Pittel,Borgs-at-al}. On the other hand, for~$\alpha>1$ one should be able to sharpen the control of the rate function near its minimum to derive a CLT for the fluctuations of the size of the giant component.
\end{remarks}

%\smallskip
Several ingredients enter our proof of Theorem~\ref{thm-main} which are of independent interest. We state these as separate theorems. The first one concerns the exponential decay rate for the probability that the random-graph is (completely) connected:

\begin{theorem}
\label{thm-connected}
Let~$K$ denote the event that $\scrG(n,\p)$ is connected. Then
\begin{equation}
\label{5.1}
P_{n,\alpha}(K)=(1-e^{-\alpha})^n\texte^{O(\log n)},
\qquad n\to\infty,
\end{equation}
where~$O(\log n)$ is bounded by a constant times $\log n$ uniformly on compact sets of~$\alpha\in[0,\infty)$.
\end{theorem}

We remark that Theorem~\ref{thm-connected} holds with~$\texte^{O(\log n)}$ replaced by~$C(\alpha)+o(1)$, see~\cite{Bender-2} for a proof. However, the requisite steps seem far in excess of the derivation in Sect.~\ref{sec3}. Furthermore, various pieces of Theorem~\ref{thm-connected} have been discovered, apparently multiple times, in~\cite{Stepanov,Knuth,Takacs,Lomonosov,Bender-at-al}; cf also the discussion following Lemma~\ref{lemma3.3}.

Next we present a result concerning the event that~$\scrG(n,\p)$ contains no cycles. Such problems have been extensively studied under the conditions where this probability is~$O(1)$, see e.g.~\cite{Bollobas}. Our theorem concerns the large-deviation properties of this event:

\begin{theorem}
\label{thm-only-trees}
Let~$L$ be the event that~$\scrG(n,\p)$ contains no cycles. Then
\begin{equation}
%\label{}
\lim_{n\to\infty}P_{n,\alpha}(L)^{1/n}
=\begin{cases}
\alpha \exp\left(-\frac{\alpha}{2} + \frac{1}{2\alpha}\right),\quad&\text{if }\alpha>1,
\\
1,\quad&\text{otherwise}.
\end{cases}
\end{equation}
\end{theorem}

Strictly speaking, this result is not needed for the proof of our main theorem; it is actually used to derive the exponential decay for the probability of the event that~$\scrG(n,\p)$ contains only ``small'' components. Surprisingly, the decay rates for these two events are exactly the same:

\begin{theorem}
\label{thm-only-small}
Let~$L$ be the event that~$\scrG(n,\p)$ contains no cycles and let~$B_r$ be the event that there are no components larger than~$r$. Then
\begin{equation}
%\label{}
\lim_{r\to\infty}\liminf_{n\to\infty}P_{n,\alpha}(B_r)^{1/n}
=
\lim_{\epsilon\downarrow0}
\limsup_{n\to\infty}P_{n,\alpha}(B_{\epsilon n})^{1/n}
=\lim_{n\to\infty}P_{n,\alpha}(L)^{1/n}.
\end{equation}
\end{theorem}

\smallskip\noindent
\emph{Update}:
In the present paper we prove Theorem~\ref{thm-only-trees} using enumeration and generating-function techniques. Recently, a probabilistic approach has been developed by which we obtain an expansion of $P_{n,\alpha}(L_n)$ to quantities of order unity. One advantage of the new approach is that it also permits the analysis of the conditional measure $P_{n,\alpha}(\cdot|L_n)$; see~\cite{BCS-trees}.
\smallskip

To finish the discussion of our results, let us give some reason for the word ``thermodynamic'' in the title. The motivation comes from an analogy with droplet formation in systems at phase transition. Such situations have been studied extensively in the context of percolation~\cite{ACC,Cerf} and Ising (and Potts) model~\cite{DKS,Ioffe-Schonmann,BCK,Bodineau,Cerf-Pisztora} under the banner of ``Wulff construction,'' see~\cite{BIV} for a review of these matters. 

One of the principal questions underlying Wulff construction is as follows: Compute the probability---and the characteristics of typical configurations carrying the event---that a given fraction of the system is in one thermodynamic state (e.g., liquid) while the rest is in another state (e.g., gas). It turns out that the typical configurations are such that the two phases separate; a droplet of one phase ``floats'' in the other phase. The requisite probability is then given by a large-deviation expression whose rate function is composed of three parts: the ``surface'' energy and entropy of the droplet, the rate function for the probability that the droplet is all in one phase, and the rate function for the probability that the complement of the droplet is in the other phase.

In the case under study, the droplet is exactly the giant component and its weight is just the probability that all vertices in the droplet are connected to each other. The ``surface'' energy is (the~log of) the probability that no vertex inside is connected to no vertex outside; the entropy is (the log of) the number of ways to choose the corresponding number of sites. The weight of the phase outside simply amounts to the probability that all remaining components are of submacroscopic scale.
When the leading-order exponential decay rate of all of these contributions is extracted using Theorems~\ref{thm-connected}--\ref{thm-only-small}, we get a quantity that only depends on the fraction of vertices taken by the droplet. The resulting expression is the one on the right-hand side of \eqref{2.6}.

\section{Everybody connected}
\label{sec3}\noindent
The goal of this section is to prove Theorem~\ref{thm-connected}. Our proof is based on showing that the probability in~\eqref{5.1} is exactly the same probability in a related, directed graph problem.

\smallskip
For a collection of vertices~$\scrV_n=\{1,\dots,n\}$ and a set of edge probabilities~$(p_{kl})_{1\le k<l\le n}$, let~$\scrG$ be the inhomogeneous undirected random graph over~$\scrV_n$. Similarly, let~$\vec\scrG$ denote the inhomogeneous \emph{directed} complete random graph with the restriction that the two possible (directed) edges between~$k$ and~$l$ occur independently, each with probability~$p_{kl}$. To keep our notation distinct from the special case~$p_{kl}=\ffrac\alpha n$ treated throughout this paper, we will write~$P$ instead of~$P_{n,\alpha}$.

\begin{definition}
A labelled directed graph $\scrG=(\scrV,\scrE)$ is said to be \emph{grounded} at vertex $v \in \scrV$ if for every $w\in\scrV$ there exists a (directed)  path from $w$ to $v$ in $\scrE$.
\end{definition}

The identification of the two problems is now stated as follows:

\begin{lemma}
\label{lemma3.2a}
Let~$K$ be the event that $\scrG$ is connected and let~$G$ be the event that $\vec{\scrG}$ is grounded at vertex~``$1$.''
Then $P(K) = P(G)$.
\end{lemma}

\begin{proofsect}{Proof}
We use induction on the total number of edges incident with vertex~``$n$.'' Indeed, if~$p_{kn}=0$ for all~$k=1,\dots,n-1$, then~$P(K)=P(G)$ because both probabilities are zero. Now let us suppose that $P(K) = P(G)$ when~$p_{\ell n}=0$ for all~$\ell=k,\dots,n-1$ and let us prove that it also for~$p_{kn}>0$. It clearly suffices to show that the partial derivatives of~$P(K)$ and~$P(G)$ with respect to~$p_{kn}$ are equal for all~$p_{kn}\in[0,1]$.

Notice first that both~$K$ and~$G$ are increasing events. Invoking Russo's formula, see~\cite{Russo} or~\cite[Theorem~2.25]{Grimmett}, we obtain
\begin{equation}
%\label{}
\frac\partial{\partial p_{kn}}P(G)=P\bigl((n,k)\text{ is pivotal for }G\bigr),
\end{equation}
where the event $\{(n,k)\text{ is pivotal for }G\}$ means that if~$(n,k)$ is occupied, the event~$G$ occurs and if not, it does not. (Note that~$(n,k)$ denotes the edge going from~``$n$'' to~``$k$.'') The conditions under which this event occurs are straightforward: The set $\scrV_n=\{1,\dots,n\}$ splits into two disjoint components, one rooted at~``$1$'' and the other at~``$n$,'' such that no vertex in the component associated with vertex~``$n$'' has an oriented edge to the other component and~$k$ has an oriented path to~$1$.
Similarly, we have
\begin{equation}
\label{3.2a}
\frac\partial{\partial p_{kn}}P(K)=P\bigl((n,k)\text{ is pivotal for }K\bigr).
\end{equation}
Here~$\{(n,k)\text{ is pivotal for }K\}$ simply means that, if the edge~$(n,k)$ is absent,~$\scrV_n$ consist of two connected components, one containing~``$1$'' and the other containing~``$n$.''

To see the equality of partial derivatives, we split both ``pivotal'' events according to the component containing the vertex~``$n$.'' If~$\scrW$ is a set of vertices such that~$n\in\scrW$ and~$1\not\in\scrW$, let~$\scrG_{n,\scrW}$ and~$\scrG_{1,\scrW}$ be the restrictions of~$\scrG$ to~$\scrW$, and $\scrV_n\setminus\scrW$, respectively. Similarly, let~$\vec\scrG_{n,\scrW}$ and~$\vec\scrG_{1,\scrW}$ be the corresponding ``components'' of the oriented graph. Let~$K_{n,\scrW}$ and $K_{1,\scrW}$ be the events that~$\scrG_{n,\scrW}$ and~$\scrG_{1,\scrW}$ are connected and let~$G_{n,\scrW}$ and $G_{1,\scrW}$ be the events that~$\vec\scrG_{n,\scrW}$ is grounded at~``$n$'' and that~$\vec\scrG_{1,\scrW}$ is grounded at~``$1$,'' respectively. Since these pairs of events are independent, we have
\begin{equation}
\label{3.3a}
P\bigl((n,k)\text{ is pivotal for }G\bigr)
=\sum_{\begin{subarray}{c}
\scrW\colon n\in\scrW\\1,k\not\in\scrW
\end{subarray}}
P(G_{1,\scrW})P(G_{n,\scrW})P(C_{\scrW})\bigl|_{p_{kn}=0},
\end{equation}
where~$C_{\scrW}$ is the event that no vertex in~$\scrW$ has a (directed) edge to~$\scrV_n\setminus\scrW$.
But the induction assumption tells us that~$P(G_{1,\scrW})=P(K_{1,\scrW})$ and~$P(G_{n,\scrW})=P(K_{n,\scrW})$, and the symmetry of edge probabilities for the directed graph tells us that~$P(C_{\scrW})$ is the probability that~$\scrG_{n,\scrW}$, and~$\scrG_{1,\scrW}$ are not connected by an edge in~$\scrG$. Substituting these into \eqref{3.3a}, we get the right-hand side of \eqref{3.2a}. This completes the induction step.
\end{proofsect}

From now on, let~$K$ and~$G$ pertain to the specific random graphs~$\scrG(n,\p)$ and~$\vec\scrG(n,\p)$. We begin with upper and lower bounds on~$P_{n,\alpha}(K)$:

\begin{lemma}
\label{lemma3.3}
$P_{n,\alpha}(K) \leq \left(1-\left(1-\p\right)^{n-1}\right)^{n-1}$.
\end{lemma}

\begin{proofsect}{Proof}
Let~$E$ be the event---concerning the graph~$\vec\scrG(n,\p)$---that every vertex except number~``$1$'' has at least one outgoing edge. Then~$G\subset E$ and so
\begin{equation}
%\label{}
P_{n,\alpha}(G) \le P_{n,\alpha}(E)=\left(1-\left(1-\p\right)^{n-1}\right)^{n-1}.
\end{equation}
Invoking Lemma~\ref{lemma3.2a}, this proves the desired upper bound. 
\end{proofsect}

We remark that the upper bound in Lemma~\ref{lemma3.3} has been discovered (and rediscovered) several times in the past. It seems to have appeared in~\cite{Stepanov} for the first time and later in~\cite{Knuth,Takacs} and also~\cite{Karp}. A generalization to arbitrary connected graphs has been achieved in~\cite{Lomonosov}.

\begin{lemma}
\label{lemma3.4}
$P_{n,\alpha}(K) \geq \left(1-\left(1-\p\right)^{n-1}\right)^{n-1} \frac1n$.
\end{lemma}

\begin{proofsect}{Proof}
Consider the following events for directed random graph~$\vec{\scrG}(n,\p)$:
Let~$E$ be the event that every vertex, except vertex number~``$1$,'' has at least one outgoing edge and let $F$ be the event every such vertex has \emph{exactly} one outgoing edge. Since~$G\subset E$, we have
\begin{equation}
%\label{}
P_{n,\alpha}(G)= P_{n,\alpha}(E) P_{n,\alpha}(G|E).
\end{equation}
We claim that
\begin{equation}
\label{3.6a}
P_{n,\alpha}(G|E)\geq P_{n,\alpha}(G|F).
\end{equation}
Indeed, let us pick an outgoing edge for each vertex different from~``$1$,'' uniformly out of all edges going out of that vertex, and let us color these edges red. Let~$G'$ be the event that~$G$ occurs using only the red edges.
The distribution of red edges conditional on~$E$ is the same as conditional on~$F$. Hence~$P_{n,\alpha}(G|E)\ge P_{n,\alpha}(G'|E)=P_{n,\alpha}(G'|F)$. But, on $F$, every available edge is red and so $P_{n,\alpha}(G'|F)=P_{n,\alpha}(G|F)$. Combining these inequalities, \eqref{3.6a} is proved.

The number of configurations that $\vec{\scrG}(n,\p)$ can take on~$F$ is exactly $(n-1)^{n-1}$.  On the other hand, the number of configurations which result in $\vec{\scrG}(n,\p)$ being grounded is $a_n = n^{n-2}$---the number of labelled trees with~$n$ vertices. 
Hence
\begin{equation}
%\label{}
P_{n,\alpha}(G|F)\ge\frac{n^{n-2}}{(n-1)^{n-1}}\ge\frac1n.
\end{equation}
Using that $P_{n,\alpha}(E)= (1-(1-\p)^{n-1})^{n-1}$ the desired bound follows.
\end{proofsect}

\begin{proofsect}{Proof of Theorem~\ref{thm-connected}}
The claim is proved by noting
\begin{equation}
\lim_{n \to \infty} \frac{\left(1-\left(1-\p\right)^{n-1}\right)^{n-1}}{\left(1-e^{-\alpha}\right)^{n-1}} = \exp \left((1-\ffrac\alpha2)\frac{\alpha e^{-\alpha}}{1-e^{-\alpha}}\right)
\end{equation}
and using the results of Lemmas~\ref{lemma3.3} and~\ref{lemma3.4}.
\end{proofsect}

\section{Only trees}
\noindent
Here we will assemble the necessary ingredients for the proof of Theorem~\ref{thm-only-trees}.
The proof is based on somewhat detailed combinatorial estimates and arguments using generating functions. 

Recall that~$L$ denotes the event that~$\scrG(n,\p)$ contains no cycles and that~$B_r$ denotes the event that all components of~$\scrG(n,\p)$ have no more than~$r$ vertices.
We begin by a combinatorial representation of the probability~$P_{n,\alpha}(L\cap B_r)$: Let~$a_\ell$ denote the number of labeled trees on~$\ell$ vertices. Then
\begin{equation}
\label{4.2}
\begin{aligned}
P_{n,\alpha}(L\cap B_r)
&=
\sum_{\begin{subarray}{c}
\sum m_\ell \ell = n  \\ m_\ell=0\,\forall\ell>r
\end{subarray}}
\frac{n!}{\prod_\ell\bigl[m_\ell ! (\ell!)^{m_\ell}\bigr]} \left(\, \prod_{\ell\ge1} \left[a_\ell \left(\frac{\alpha}{n}\right)^{\ell-1}\right]^{m_\ell}\right) \left(1 - \frac{\alpha}{n}\right)^{\binom{n}{2} - n + \sum m_\ell}
\\
&=n! \left(\frac{\alpha}{n}\right)^n \left(1 - \frac{\alpha}{n}\right)^{\binom{n}{2} - n} \sum_{k=1}^n \left(\frac{\alpha}{n}\right)^{-k} \left(1 - \frac{\alpha}{n}\right)^k 
Q_{n,k,r},
\end{aligned}
\end{equation}
where we set~$k=\sum_\ell m_\ell$, applied the constraint~$\sum_\ell \ell m_\ell=n$ and let~$Q_{n,k,r}$ denote the sum
\begin{equation}
%\label{}
Q_{n,k,r} 
=\sum_{\substack{\sum m_\ell \ell = n \\ \sum m_\ell = k \\ m_\ell = 0\,\forall \ell > r}} \prod_{\ell\ge1}\left( \frac{a_\ell}{\ell!}\right)^{m_\ell} \frac{1}{m_\ell !}.
\end{equation}
We begin by isolating the large-$n,k$ behavior of this quantity:

\begin{proposition}
%\label{lemma}
Consider the polynomial
\begin{equation}
%\label{}
F_r(s) = \sum_{\ell=1}^r \frac{s ^\ell a_\ell}{\ell!}
\end{equation}
Then for all~$n,k,r\ge1$,
\begin{equation}
%\label{}
\label{Q-F}
Q_{n,k,r}\le \frac{1}{k!} \inf_{s>0} \frac{F_r(s)^k}{s^n}.
\end{equation}
Moreover, for each~$\eta>0$, there is~$n_0<\infty$ and a sequence~$(c_r)_{r\ge1}$ of positive numbers for which
\begin{equation}
\label{Q-F2}
Q_{n,k,r}\ge\frac{c_r}{\sqrt n}\frac{1}{k!} \inf_{s>0} \frac{F_r(s)^k}{s^n}
\end{equation}
holds for all~$n\ge n_0$, all~$k\ge1$ and all~$r\ge2$ such that~$k<(1-\eta)n$ and~$rk>n(1+\eta)$.
\end{proposition}

\begin{proofsect}{Proof of upper bound}
Let us consider the generating function
\begin{equation}
%\label{}
\hat{Q}_r(s, z) = 1 + \sum_{n=1}^\infty \sum_{k=1}^n Q_{n,k,r} z^k s^n
=\exp \bigl\{z F_r(s)\bigr\},
\end{equation}
where we used Fubini-Tonelli to derive the second equality. 
Since~$F_r$ is a polynomial, the Cauchy integral formula yields
\begin{equation}
\label{4.6a}
Q_{n,k,r} =\frac{1}{(2\pi\texti)^2} \oint \textd s \oint\textd z \frac{\exp\{z F_r(s)\}}{s^{n+1}z^{k+1}}
=\frac{1}{2\pi\texti}\frac1{k!} \oint \textd s\frac{F_r(s)^k}{s^{n+1}},
\end{equation}
where all integrals are over a circle of positive radius centered at the origin of~$\C$.
Since all coefficients of~$F_r$ are non-negative, $\theta\mapsto |F_r(s\texte^{\texti\theta})|$ for~$s>0$ is maximized at~$\theta=0$. Bounding the integrand by its value at~$\theta=0$, the integral yields a factor~$2\pi$; optimizing over~$s>0$ then gives the upper bound in \eqref{Q-F}.
\end{proofsect}

\begin{proofsect}{Proof of lower bound}
As is common in Tauberian arguments, the lower bound will require somewhat more effort. First let us note that under the conditions $k<(1-\eta)n$ and~$rk>n(1+\eta)$ the function~$s\mapsto F_r(s)^k/s^n$, for~$s>0$, blows up both at~$0$ and~$\infty$. Its minimum is thus achieved at an interior point; for the rest of this proof we will fix~$s$ to a minimizer of this function. Since~$|F_r(s\texte^{\texti\theta})|<F_r(s)$ for all~$\theta\in(-\pi,\pi]\setminus\{0\}$, the part of the integral in \eqref{4.6a} corresponding to~$|\theta|>\epsilon$ is exponentially small (in~$n$) compared to the infimum in~\eqref{Q-F2}. We thus need to show the lower bound only for the portion of the integral over~$\theta$ with~$|\theta|\le\epsilon$, for some fixed~$\epsilon>0$. 

Since~$F_r$ has positive coefficients,~$F_r\ne0$ in the (complex) $\epsilon$-neighborhood of~$s$. This allows us to define the function
\begin{equation}
g(\theta)=\log\frac{F_r(s\texte^{\texti\theta})^\varrho}{s\texte^{\texti\theta}},\qquad|\theta|\le\epsilon,
\end{equation}
where~$\varrho$ plays the role of~$\ffrac kn$. The function~$g$ is analytic in an $O(\epsilon)$-neighborhood of the origin. The choice of~$s$ implies that~$g'(0)=0$ which is equivalent to
\begin{equation}
%\label{}
\frac{sF_r'(s)}{F_r(s)}=\frac1\varrho.
\end{equation}
For the second derivative we get $g''(0)=-\varrho\text{Var}(X)$, where~$X$ is the random variable with law
\begin{equation}
%\label{}
P(X=\ell)=\frac1{F_r(s)}\frac{a_\ell s^\ell}{\ell!},\qquad \ell=1,\dots,r.
\end{equation}
In particular, since our restrictions on~$\varrho$ between~$\frac1r(1+\eta)$ and~$1-\eta$ imply that~$s$ is bounded away from zero, this law is non-degenrate and so~$g''(0)<0$.

The analyticity of~$\theta\mapsto g(\theta)$ for~$\theta=O(\epsilon)$ implies that~$g'''$ is bounded in this neighborhood, and so by Taylor's theorem we have
\begin{equation}
%\label{}
g(\theta)=g(0)-A\theta^2+O(\theta^3),
\end{equation}
where~$A=A(r,\varrho)$ is positive uniformly in the allowed range of~$\varrho$'s and~$O(\theta^3)$ is a quantity bounded by~$|\theta|^3$ times a constant depending only on~$r$, $\epsilon$ and~$\eta$. (In particular, we may assume that $O(\theta^3)$ is dominated by~$\frac12A\theta^2$ for~$|\theta|\le\epsilon$.)

We will split the integral over~$\theta\in[-\epsilon,\epsilon]$ into two more parts. Let~$\delta>0$ and note that~$ng(0)$ is the logarithm of the infimum in~\eqref{Q-F2}. Then for~$\theta$ with~$\delta n^{-1/3}\le |\theta|\le\epsilon$ we have
\begin{equation}
%\label{}
n\RE g(\theta)\le n g(0)-\frac12A\delta^2 n^{1/3}
\end{equation}
which shows that even this portion of the integral brings a contribution that is negligible compared to the right-hand side of~\eqref{Q-F2}.
But for~$|\theta|\le\delta n^{-1/3}$ we have $nO(\theta^3)=O(\delta)$ and so for~$\delta\ll1$, the Taylor remainder will always have imaginary part between, say,~$-\ffrac\pi4$ and~$\ffrac\pi4$. This means that
\begin{equation}
%\label{}
\RE\int_{-\delta n^{-1/3}}^{\delta n^{-1/3}} \texte^{ng(\theta)}\,\textd\theta
\ge\frac12\texte^{ng(0)}\int_{-\delta n^{-1/3}}^{\delta n^{-1/3}}\texte^{-nA\theta^2}\,\textd\theta
\ge \frac{c}{\sqrt n}\,\texte^{ng(0)}
\end{equation}
for some constant~$c>0$ which may depend on~$r$ and~$\eta$ but not on~$\varrho$ and~$n$. Combined with the previous estimates, this proves the lower bound~\eqref{Q-F2}.
\end{proofsect}

In light of the above lemma, the $k$-th term in the sum on the extreme right of \eqref{4.2} becomes
\begin{equation}
\label{4.14a}
\alpha^{-k} n^k\texte^{-\alpha\frac{k}{n}} Q_{n,k,r}
=\texte^{o(n)}\inf_{s>0}\exp\bigl\{n\Theta_r(s,\ffrac kn)\bigr\},
\end{equation}
where
\begin{equation}
%\label{}
\Theta_r(s, \varrho)=-\varrho \log \alpha - \varrho \log \varrho + \varrho + \varrho \log F_r(s) - \log s.
\end{equation}
Here we should interpret \eqref{4.14a} as an upper bound for~$r=n$ and a lower bound for fixed~$r$. It is clear that, regardless of~$r$, the sum is dominated by $k=\lfloor\varrho n\rfloor$ for which $\varrho\mapsto\inf_{s>0}\Theta_r(s,\varrho)$ is maximal. Such values are characterized as follows:

\begin{lemma}
\label{lemma4.2}
Let~$\alpha>0$ and~$r\ge2$. Then there is a unique~$(s_r,\varrho_r)\in[0,\infty]\times[\ffrac1r,1]$ for which
\begin{equation}
\label{4.16a}
\Theta_r(s_r,\varrho_r)=\sup_{1/r\le\varrho\le1}\inf_{s>0}\,\Theta_r(s,\varrho).
\end{equation}
Moreover, we always have~$s_r\in(0,\infty)$ and~$\varrho_r\in(\ffrac1r,1)$ and, furthermore,
\begin{equation}
\label{4.17a}
\lim_{r\to\infty}\Theta_r(s_r,\varrho_r)=\begin{cases}
1+\ffrac\alpha2-\log\alpha,\qquad&\text{if }\alpha\le1,
\\*[1mm]
1+\frac1{2\alpha},\qquad&\text{if }\alpha>1.
\end{cases}
\end{equation}
\end{lemma}

\begin{proofsect}{Proof}
We begin by ruling out the ``boundary values'' of~$s$ and~$\varrho$. First, if~$\varrho=\ffrac1r$, then the infimum over~$s$ is actually achieved by~$s=\infty$. In that case~$F_r(s)=\infty$ and the (one-sided) derivative with respect to~$\varrho$ is infinite, i.e.,~$\varrho=\ffrac1r$ is a strict local minimum of $\varrho\mapsto\inf_{s>0}\Theta_r(s,\varrho)$. Similarly, for~$\varrho=1$ the infimum over~$s>0$ is achieved at~$s=0$ but then the $\varrho$-derivative of $\varrho\mapsto\inf_{s>0}\Theta_r(s,\varrho)$ is negative infinity, i.e., also~$\varrho=1$ is a strict local minimum. It follows that any $(s_r,\varrho_r)$ satisfying \eqref{4.16a} necessarily lies in $(0,\infty)\times(\ffrac1r,1)$.

Setting the partial derivatives with respect to~$s$ and~$\varrho$ to zero shows that any minimizing pair is the solution of the equations
\begin{equation}
F_r(s) = \alpha \varrho 
\quad\text{and}\quad
s F_r'(s) = \alpha.
\end{equation}
In light of monotonicity of~$s\mapsto sF_r'(s)$, the solution is actually unique. To figure out the asymptotic as~$r\to\infty$, we note that for~$s\le\ffrac1\texte$, 
\begin{equation}
%\label{}
sF_r'(s) = \sum_{\ell = 1}^{r} a_{\ell}\frac{s^{\ell}}{(\ell-1)!}\,\underset{r\to\infty}\longrightarrow\, W(s),
\end{equation}
where~$W$ is the unique number in~$[0,\ffrac1\texte]$ such that~$W\texte^{-W}=s$. (Incidentally, $W$ is closely related to the survival probability of the Galton-Watson branching process with Poisson offspring distribution.) If~$s>\ffrac1\texte$, then~$sF_r'(s)\to\infty$ as~$r\to\infty$. 
Using the relation between~$sF'_r(s)$ and~$\alpha$, we thus get
\begin{equation}
s_r\,\underset{r\to\infty}\longrightarrow\,\begin{cases}
\alpha\texte^{-\alpha},\qquad&\text{if }\alpha\le1,
\\*[1mm]
\ffrac1\texte,\qquad&\text{if }\alpha>1.
\end{cases}
\end{equation}
Integrating the derivative of~$F_r$ now shows that $F_r(s_r)\to\alpha (1-\ffrac{\alpha}{2})$ for $\alpha \le 1$. Using that $F_r'(s)$ is bounded for~$s\le s_r$, we also find that $F_r(s_r)\to\ffrac12$ for $\alpha \ge1$. This yields
\begin{equation}
\varrho_r\,\underset{r\to\infty}\longrightarrow\,\begin{cases}
1-\ffrac\alpha2,\qquad&\text{if }\alpha\le1,
\\*[1mm]
\frac1{2\alpha},\qquad&\text{if }\alpha>1.
\end{cases}
\end{equation}
Noting that~$\Theta(s_r,\varrho_r)=\varrho_r-\log s_r$ we now get \eqref{4.17a}.
\end{proofsect}

\begin{proofsect}{Proof of Theorem~\ref{thm-only-trees}}
By the fact that the supremum over~$\varrho$ in \eqref{4.16a} is achieved at an interior point, we can control the difference between the maximizing~$\ffrac kn$ and its continuous counterpart~$\varrho$. Thence
\begin{equation}
%\label{}
P_{n,\alpha}(L\cap B_r)=q_{n,r}
\, n! \left(\frac{\alpha}{n}\right)^n 
\texte^{-\alpha n/2} 
\,\exp\bigl\{n\Theta_r(s_r, \varrho_r)\bigr\},
\end{equation}
where
\begin{equation}
%\label{}
\frac{\tilde c_r}{\sqrt n}\le q_{n,r}\le n
\end{equation}
for some positive constants~$\tilde c_r$ which may depend on~$r$ and~$\alpha$.
Since~$B_n$ contains every realization of~$\scrG(n,\p)$, taking~$r=n$ and applying Lemma~\ref{lemma4.2} directly shows that~$P_{\alpha,n}(L)\le\texte^{n\Psi(\alpha)+o(n)}$. To get a corresponding lower bound, we fix~$r\ge2$ and apply~$P_{\alpha,n}(L)\ge P_{\alpha,n}(L\cap B_r)$. Taking~$\ffrac1n$-th power and letting~$n\to\infty$ then yields
\begin{equation}
%\label{}
\lim_{n\to\infty}P_{n,\alpha}(L\cap B_r)^{1/n}
=\alpha\texte^{-1-\alpha/2+\Theta_r(s_r,\varrho_r)}.
\end{equation}
As we have just checked, the right-hand side tends to~$\texte^{\Psi(\alpha)}$ as~$r\to\infty$.
\end{proofsect}

\begin{corollary}
\label{cor-BLL}
We have
\begin{equation}
%\label{}
\lim_{r\to\infty}\lim_{n\to\infty}P_{\alpha,n}(B_r\cap L)^{1/n}=\lim_{n\to\infty}P_{\alpha,n}(L)^{1/n}.
\end{equation}
\end{corollary}

\begin{proofsect}{Proof}
This summarizes the last step of the previous proof.
\end{proofsect}

\section{No big = no cycles}
\noindent
Here we will prove that absence of large component has a comparable cost to absence of cycles, at least on an exponential scale. To achieve this goal, apart from Corollary~\ref{cor-BLL}, we will need the following upper bound:

\begin{lemma}
\label{lemma5.1}
Let $B_r$ be the event that~$\scrG(n,\ffrac\alpha n)$ has no components larger than~$r$ and let~$L$ be the event that all connected components of~$\scrG(n,\ffrac\alpha n)$ are trees. Then for all~$r\ge1$,
\begin{equation}
\label{NB-NL}
P_{n,\alpha}(B_r)\le P_{n,\alpha}(L)\Bigl(1-\frac\alpha n\Bigr)^{-\frac12rn}.
\end{equation}
\end{lemma}

\begin{proofsect}{Proof}
Let~$C$ be the restriction of~$\scrG(n,\p)$ to a set~$S\subset\{1,\dots,n\}$. Let~$T$ be a tree on~$S$. Then
\begin{equation}
%\label{}
\frac{P_{n,\alpha}(C=T)}{P_{n,\alpha}(C\supset T)}=\Bigl(1-\frac\alpha n\Bigr)^{\binom{|S|}2-|S|+1}
\ge\Bigl(1-\frac\alpha n\Bigr)^{\frac12|S|^2}.
\end{equation}
Hence
\begin{equation}
\label{5.3a}
P_{n,\alpha}(C\text{ is connected})\le\sum_TP_{n,\alpha}(C\supset T)\le\Bigl(1-\frac\alpha n\Bigr)^{-\frac12|S|^2}
P_{n,\alpha}(C\text{ is a tree}). 
\end{equation}
Now, if~$L_r$ is the event that no component of~$\scrG(n,\p)$ of size larger than~$r$ has cycles, then $B_r\subset L_r$ and so $P_{n,\alpha}(B_r)\le P_{n,\alpha}(L_r)$. Let~$\{S_j\}$ be a partition of~$\{1,\dots,n\}$ and let~$P_{n,\alpha}(\{S_j\})$ denote the probability that~$\{S_j\}$ are the connected components of~$\scrG(n,\p)$. Then
\begin{equation}
\label{5.4a}
P_{n,\alpha}(L_r)=\sum_{\{S_j\}}P_{n,\alpha}(\{S_j\})P_{n,\alpha}(L_r|\{S_j\}),
\end{equation}
where~$P_{n,\alpha}(L_r|\{S_j\})$ is the conditional probability of~$L_r$ given that~$\{S_j\}$ are the connected components of~$\scrG(n,\p)$.

Letting~$C_j$ represent the restriction of~$\scrG(n,\p)$ to~$S_j$, the bound \eqref{5.3a} tells us that
\begin{equation}
\begin{aligned}
P_{n,\alpha}(L_r|\{S_j\})
&=\prod_{j\colon |S_j|\ge r}P_{n,\alpha}(C_j\text{ is a tree}|C_j\text{ is connected})
\\
&\le\,\prod_jP_{n,\alpha}(C_j\text{ is a tree}|C_j\text{ is connected})\prod_{j\colon |S_j|<r}
\Bigl(1-\frac\alpha n\Bigr)^{-\frac12|S_j|^2}
\end{aligned}
\end{equation}
Using that~$|S_j|<r$ for every~$S_j$ contributing to the second product and applying that the sum of~$|S_j|$ over the components with $|S_j|<r$ gives at most~$n$, we then get
\begin{equation}
%\label{}
P_{n,\alpha}(L_r|\{S_j\})\le P_{n,\alpha}(L|\{S_j\})\Bigl(1-\frac\alpha n\Bigr)^{-\frac12rn}.
\end{equation}
Plugging this back in \eqref{5.4a}, the desired bound follows.
\end{proofsect}

\begin{proofsect}{Proof of Theorem~\ref{thm-only-small}}
By Lemma~\ref{lemma5.1} we have
\begin{equation}
%\label{}
%\lim_{n\to\infty}P_{n,\alpha}(B_r)^{1/n}
\limsup_{n\to\infty}P_{n,\alpha}(B_{\epsilon n})^{1/n}
\le \texte^{\epsilon/2}
\lim_{n\to\infty}P_{n,\alpha}(L)^{1/n}.
\end{equation}
On the other hand, the inclusion~$B_r\supset B_r\cap L$ and Corollary~\ref{cor-BLL} yield
\begin{equation}
%\label{}
\liminf_{n\to\infty} P_{n,\alpha}(B_r)^{1/n}\ge 
\liminf_{n\to\infty} P_{n,\alpha}(B_r\cap L)^{1/n}\,\underset{r\to\infty}\longrightarrow\, 
\lim_{n\to\infty}P_{n,\alpha}(L).
\end{equation}
Since~$P_{n,\alpha}(B_r)\le P_{n,\alpha}(B_{\epsilon n})$ eventually for any fixed~$r\ge1$ and~$\epsilon>0$, all limiting quantities are equal provided we take~$r\to\infty$ and/or~$\epsilon\downarrow0$ after~$n\to\infty$.
\end{proofsect}

\section{Proof of main result}
\noindent
Before we start proving our main result, Theorem~\ref{thm-main}, we need to ensure that if a large component is present in the graph, then it is unique. The statement we need is as follows:

\begin{lemma}
\label{lemma6.1}
Let~$K_{\epsilon,2}$ be the event that~$\scrG(n,\p)$ is either connected or has exactly two connected components, each of which is of size at least~$\epsilon n$, and recall that~$K$ is the event that~$\scrG(n,\p)$ is connected. Then for all~$\alpha_0>0$ and~$\epsilon_0>0$ there exists~$c_1=c_1(\alpha_0,\epsilon_0)<1$ such that for all~$\epsilon\ge\epsilon_0$ and all~$\alpha\le\alpha_0$,
\begin{equation}
%\label{}
\limsup_{n\to\infty}P_{\alpha,n}(K^\cc|K_{\epsilon,2})^{1/n}<c_1.
\end{equation}
\end{lemma}

\begin{proofsect}{Proof}
It clearly suffices to show that the ratio of~$P_{\alpha,n}(K_{\epsilon,2}\setminus K)$ and~$P_{\alpha,n}(K)$ decays to zero exponentially with~$n$, with a rate that is uniformly bounded in~$\epsilon\ge\epsilon_0$ and~$\alpha\le\alpha_0$. In light of Theorem~\ref{thm-connected} and the fact that~$K_{\epsilon,2}$ only admits components that grow linearly with~$n$, we have
\begin{equation}
\label{6.3}
\frac{P_{\alpha,n}(K_{\epsilon,2}\setminus K)}{P_{\alpha,n}(K)}
=\,\texte^{o(n)}\!\!
\sum_{\epsilon n\le k\le n-\epsilon n}
\binom nk\,\frac{\pi_1(\alpha\,\ffrac kn)^k\,\pi_1\bigl(\alpha(1-\ffrac kn)\bigr)^{n-k}}{\pi_1(\alpha)^n}\Bigl(1-\frac\alpha n\Bigr)^{k(n-k)},
\end{equation}
where~$o(n)/n$ tends to zero uniformly in~$\epsilon\ge\epsilon_0$ and~$\alpha\le\alpha_0$.
Writing~$\varrho$ for~$\ffrac kn$, the expression under the sum can be bounded by~$\texte^{n[\Xi(\varrho)-\Xi(0)]}$, where
\begin{equation}
%\label{}
\Xi(\varrho)=-S(\varrho)+\varrho\log\pi_1(\alpha\varrho)+(1-\varrho)\log\pi_1\bigl(\alpha(1-\varrho)\bigr)-\alpha\varrho(1-\varrho).
\end{equation}
Since~$\varrho$ is restricted to the interval~$[\epsilon,1-\epsilon]$, the right-hand side of \eqref{6.3} will be exponentially small if we can show~$\Xi(\varrho)<\Xi(0)$ for all~$\varrho\in(0,1)$ and all~$\alpha$.

As is easy to check, the function~$\varrho\mapsto\Xi(\varrho)$ is symmetric about the point~$\varrho=\ffrac12$. Hence, if we can prove that it is strictly convex throughout~$[0,1]$, then it is maximized at the endpoints. Introducing the function
\begin{equation}
%\label{}
G(\eta)=\eta\log\frac{\pi_1(\eta)}\eta
\end{equation}
we have
\begin{equation}
%\label{}
\alpha\Xi(\ffrac\eta\alpha)=G(\eta)+G(\alpha-\eta)+\eta(\alpha-\eta).
\end{equation}
In order to prove strict convexity of~$\Xi$, it thus suffices to show that~$G''(\eta)+1>0$ for all~$\eta>0$. Introducing yet another abbreviation~$q(\eta)=\eta/(1-\texte^{-\eta})$, a tedious but straightforward differentiation yields
\begin{equation}
%\label{}
G''(\eta)+1=\frac1q(q'-q)(q\texte^{-\eta}-1).
\end{equation}
A direct evaluation now shows that both~$q'-q$ and~$q\texte^{-\eta}-1$ are negative once~$\eta>0$.
\end{proofsect}

We will use the above lemma via the following simple conclusion:

\begin{lemma}
\label{cor6.2}
Let~$N_r$ denote the number of connected components of~$\scrG(n,\p)$ of size at least~$r$ and let~$\scrV_r$ be the set of vertices contained in these components. Then for all $\alpha>0$ and $\varrho>\epsilon>0$ there exists $c=c(\epsilon,\varrho,\alpha)>0$ such that 
\begin{equation}
\label{6.1}
P_{\alpha,n}\bigl(|\scrV_{\epsilon n}|=\lfloor\varrho n\rfloor \,\&\, N_{\epsilon n}=1\bigr)
\ge (1-\texte^{-cn})P_{\alpha,n}\bigl(|\scrV_{\epsilon n}|=\lfloor\varrho n\rfloor\bigr).
\end{equation}
\end{lemma}

\begin{proofsect}{Proof}
Clearly, \eqref{6.1} will follow if we can prove that
\begin{equation}
\label{6.8}
P_{\alpha,n}\bigl(|\scrV_{\epsilon n}|=\lfloor\varrho n\rfloor \,\&\, N_{\epsilon n}>1\bigr)
\le \texte^{-cn} P_{\alpha,n}\bigl(|\scrV_{\epsilon n}|=\lfloor\varrho n\rfloor\bigr).
\end{equation}
Let~$\scrV(x)$ denote the connected component of~$\scrG(n,\p)$ containing the vertex~$x$ and let~$x\nleftrightarrow y$ denote the event that~$x,y\in\scrV_{\epsilon n}$ but $\scrV(x)\cap\scrV(y)=\emptyset$.  Then \eqref{6.8} will be proved once we show
\begin{equation}
\label{6.9a}
P_{\alpha,n}\bigl(|\scrV_{\epsilon n}|=\lfloor\varrho n\rfloor \,\&\, x\nleftrightarrow y\bigr)
\le \texte^{-2cn}P_{\alpha,n}\bigl(|\scrV_{\epsilon n}|=\lfloor\varrho n\rfloor\bigr).
\end{equation}
(Indeed, the sum over~$x,y$ adds only a multiplicative factor of~$n^2$ on the right-hand side.)
By conditioning on the set~$\scrV_{\epsilon n}$ \emph{and} the set~$\scrV(x)\cup\scrV(y)$, this inequality will in turn follow from
\begin{equation}
\label{6.10}
P_{\alpha,n}\bigl(x\nleftrightarrow y\,\&\, \scrV(x)\cup\scrV(y)=\scrV\bigr)
\le \texte^{-2cn}P_{\alpha,n}\bigl(\scrV(x)\cup\scrV(y)=\scrV\bigr).
\end{equation}
Indeed, let us multiply both sides by the probability that~$\scrV$ is disconnected from the rest of the graph and that all components disjoint from~$\scrV$ of size at least~$\epsilon n$ take the total volume~$\lfloor\varrho n\rfloor-|\scrV|$. The sum over all admissible~$\scrV$ reduces \eqref{6.10} to \eqref{6.9a}.

We will deduce \eqref{6.10} from Lemma~\ref{lemma6.1}. Recall that~$K$ is the event that the graph is connected and~$K_{\epsilon,2}$ is the event that it has at most two components, each of which is of size at least~$\epsilon n$. We will now use these events for the restriction of~$\scrG(n,\p)$ to~$\scrV$: Let~$m=|\scrV|$, $\tilde\alpha=\alpha\frac mn$ and~$\tilde\epsilon=\epsilon\frac nm$. Then we have
\begin{equation}
%\label{}
\{x\nleftrightarrow y\}\cap\bigl\{\scrV(x)\cup\scrV(y)=\scrV\bigr\}\subset K^\cc\cap K_{\tilde\epsilon,2},
\end{equation}
while for the event on the right-hand side of \eqref{6.10} we simply get
\begin{equation}
%\label{}
\bigl\{\scrV(x)\cup\scrV(y)=\scrV\bigr\}=K_{\tilde\epsilon,2}.
\end{equation}
By Lemma~\ref{lemma6.1} and the fact that~$\tilde\alpha\le\alpha$ and $\tilde\epsilon\ge\epsilon$, 
\begin{equation}
%\label{}
P_{\tilde\alpha,m}(K^\cc|K_{\tilde\epsilon,2})\le\texte^{-c_1m},
\end{equation}
once~$n$ is sufficiently large. 
But~$m\ge2\epsilon n$ and so \eqref{6.10} holds with~$c=\epsilon c_1$.
\end{proofsect}

Now we have finally amassed all ingredients needed for the proof of our main result.

\begin{proofsect}{Proof of Theorem~\ref{thm-main}}
The case~$\varrho=0$ is quickly reduced to Theorems~\ref{thm-only-trees}--\ref{thm-only-small} while~$\varrho=1$ boils down to Theorem~\ref{thm-connected}. Thus, we are down to the cases~$\varrho\in(0,1)$. Let~$\epsilon\in(0,\varrho)$. By Lemma~\ref{cor6.2}, we can focus on the situations with~$N_{\epsilon n}=1$. To make our notation simple, let us assume that~$\varrho n$ is an integer. Then we have
\begin{equation}
%\label{}
P_{\alpha,n}\bigl(|\scrV_{\epsilon n}|=\varrho n \,\&\, N_{\epsilon n}=1\bigr)
=\binom n{\varrho n}P_{\varrho n,\alpha\varrho}(K)P_{n-\varrho n,\alpha(1-\varrho)}(B_{\epsilon n})\Bigl(1-\frac\alpha n\Bigr)^{\varrho n(1-\varrho)n}.
\end{equation}
The terms on the right-hand side represent the following: the number of ways to choose the unique component of size~$\varrho n$, the probability that this component is connected, the probability that the complement contains no component of size larger than~$\epsilon n$ and, finally, the probability that the two parts of the graph do not have any edge between them.
Invoking Stirling's formula to deal with the binomial term, and plugging explicit expressions for $P_{\varrho n,\alpha\varrho}(K)$ and $P_{n-\varrho n,\alpha(1-\varrho)}(B_{\epsilon n})$ from Theorems~\ref{thm-connected}--\ref{thm-only-small}, the result reduces to a simple calculation.
\end{proofsect}

\section*{Acknowledgments}
\noindent 
This research was partially supported by the NSF grants DMS-0306167, DMS-0301795 and DMS-0505356. We wish to thank anonymous referees for advice on style and literature.

\end{document}